\documentclass[12pt,a4paper]{amsart}

\usepackage{hyperref, color}

\newtheorem{theorem}{Theorem}[section]
\newtheorem{lemma}[theorem]{Lemma}

\theoremstyle{remark}
\newtheorem{remark}[theorem]{Remark}

\numberwithin{equation}{section}

\newcommand{\field}[1]{\mathbb{#1}}
\newcommand{\real}{\field{R}}
\newcommand{\z}{\langle}
\newcommand{\x}{\rangle}
\DeclareMathOperator{\ind}{Ind}
\DeclareMathOperator{\indw}{Ind_w}
\DeclareMathOperator{\trace}{tr}
\DeclareMathOperator{\ric}{Ric}
\DeclareMathOperator{\diver}{div}

\begin{document}

\title[The index of compact constant mean curvature surfaces]
{Lower bounds  for the index of compact  constant mean curvature surfaces in $\mathbb R^{3}$ and $\mathbb S^{3}$}
\author[Cavalcante]{Marcos P. Cavalcante}
\address{ \newline 
Instituto de Matem\'atica 
\newline Universidade Federal de Alagoas (UFAL)
\newline Campus A. C. Sim\~oes, BR 104 - Norte, Km 97, 57072-970.
\newline Macei\'o - AL -Brazil}
\email{marcos@pos.mat.ufal.br}
\thanks{The first author was partially  supported  by CNPq-Brazil.}

\author[de Oliveira]{Darlan F. de  Oliveira}
\address{
\newline Departamento de Ci\^encias Exatas 
\newline  Universidade Estadual de Feira de Santana  (UEFS)
\newline Avenida Transnordestina, S/N, Novo Horizonte, 44036-900
 \newline Feira de Santana - BA - Brazil}
\email{darlanfdeoliveira@gmail.com}


\date{\today}



\begin{abstract}
{Let $M$ be a compact constant mean curvature surface either in $\mathbb{S}^3$ or $\mathbb{R}^3$. 
In this paper we prove that the stability index of $M$ is bounded from below by a linear function of the genus. 
As a by-product we obtain a comparison theorem between the spectrum of the Jacobi operator of $M$ 
and those of Hodge Laplacian of $1$-forms on $M$..
}
\end{abstract}

\maketitle


\section{Introduction}

Let $\overline M^{3}$ be a complete Riemannian 3-manifold and let $M\subset \overline M^{3} $ be a 
compact  surface immersed in $\bar M^3$.  
It is well known that $M$ is a \emph{minimal surface} if it is  a critical point of the area functional, 
and $M$ is a \emph{constant non zero mean curvature (CMC) surface}  
if it is a critical point of the  area functional for those variations that preserve the enclosed volume. 
When the ambient space is the Euclidean three-space, it means that minimal surfaces are the mathematical models of soap films while 
constant mean curvature surfaces are those of soap bubbles. 

The stability properties of minimal and CMC surfaces are given by the study of the second variation of the area 
functional. In order to give a precise definition let us assume that $M$ is closed and two-sided, 
and denote by $N$ 
a unit normal vector field along $M$.
Given $u\in C^\infty(M)$, a smooth function on $M$ and considering variations given by $V:=uN$, 
the second variation formula is given by the quadratic form
 \[
 Q(u,u):=\int_M\|\nabla u\|^2 -(\overline {\ric}(N)+\|A\|^2)u^2 \,dM, 
 \]
where  $\overline {\ric}(N)$ is the Ricci curvature of $\overline M$ in the direction of 
$N$ and $\|A\|^2$ stands for the square norm of the second fundamental form of $M$
(see Section \ref{preliminaries} for details).
For CMC sufaces the requirement that variations preserve the enclosed volume is equivalent to
consider  functions satisfying $\int_Mu\,dM=0$. For future reference we will denote by $\mathcal F$ 
the  space of  smooth functions satisfying this property.

The \emph{Morse index} of a minimal surface $M$ is  
defined as the maximal dimension of a  vector subspace of
$C^\infty(M)$ where $Q$ is negative defined and will be  denoted by $\ind(M)$.
If $M$ is a CMC surface, we define the \emph{weak index} $\indw(M)$ as
the maximal dimension  of a vector subspace of $\mathcal F$ where $Q$ is negative defined. 
The index is always finite when $M$ is compact and if the index is zero, we say that
the surface is \emph{stable}. 
In fact, the index indicates the number of directions whose variations decrease area. 
If $M$ is non compact, we can extend the notion of
index by taking a limit of the indices of an exhaustion on $M$.
Finally, we recall that all these concepts can be given in higher dimensions as well.

A classical theorem proved independently by do Carmo and Peng \cite{dCP}, Fischer-Colbrie and Schoen
\cite{FCS} and Pogorelov \cite{P} asserts that 
\emph{a complete stable minimal surface in $\real^3$ is a flat plane}. 

If $M$ is a complete minimal surface with finite index immersed in a manifold $\overline M$ with 
nonnegative scalar curvature,  
Fischer-Colbrie \cite{FC} and Gulliver \cite{G}, also independently, proved that
$M$ is conformally diffeomorphic to a compact Riemann surface with genus
$g$  and punctured at finitely many points $p_1,\ldots, p_r$, corresponding to the ends of $M$. 
In \cite{R},  using harmonic $1$-forms to construct test functions, 
Ros proved that the index of a minimal surface
immersed in $\real^3$ or in quotients of $\real^3$ is bounded from below by a linear function of its
genus, namely he proved that $\ind(M)\geq 2g/3 $ if $M$ is orientable and  $\ind(M)\geq g/3 $ if $M$ is
nonorientable. In the case of oriented minimal surfaces  
$M\subset \real^3$, Chodosh and M\'aximo \cite{CM}, improved Ros' ideas 
and proved that $Ind(M)\geq \frac 2 3 (g+r) - 1$.

For closed minimal hypersurfaces in the  unit sphere $M\subset\mathbb S^{n+1}$,
Savo \cite{Sa} proved that the Morse index is bounded
from below by a linear function of its first Betti number (the genus in dimension 2) of $M$.
This result was  generalized recently by Ambrozio, Carlotto and Sharp in \cite{ACS} and by
Mendes and Radesh in \cite{MR} for  a wide class of ambient manifolds with positive curvature. 
These results provide a partial answer for a conjecture of Marques
and Neves which asserts that \emph{the index of a compact minimal hypersurface immersed in 
ambient spaces with positive Ricci curvature is bounded from below by a linear function
of the first Betti number} (see \cite{Ma} and \cite{N}). 

Also recently,  Chao Li   proved in \cite{Li} that the index and the nullity 
(the dimension of the subspace of solutions to $J=0$)
of a complete minimal hypersurface with  finite total curvature in the Euclidean space
is bounded from below by a linear function of the number of ends and its first Betti number.

In the case of CMC surfaces, Barbosa, Do Carmo and Eschenburg \cite{BdCE}  proved that 
\emph{geodesic spheres are the only compact constant mean curvature hypersurfaces in space forms 
that are stables}, that is, whose  weak index is zero.  
On the other hand,  there are few results  in the literature about index estimates
of other examples of CMC hypersufaces. To the best of the authors' knowledge, estimates were given 
by Lima-Sousa Neto-Rossman \cite{LSR} and Rossman \cite{Ro} for CMC tori in $\real^3$, 
by  Perdomo-Brasil  \cite{PB} for CMC hypersurfaces in the sphere, but in terms  of the dimension, 
and by Rossman-Sultana \cite{RS} and Ca\~nete \cite{Ca} for CMC tori in $\mathbb S ^{3}.$

The purpose of the present paper is to generalize Ros and Savos's estimates 
to obtain lower bounds  for the week stability index of compact CMC  surfaces immersed in either 
$\mathbb R^{3}$ or in $\mathbb S^{3}$ in terms of the  genus.
To do that we apply Ros' ideas \cite{R} making use of harmonic 1-forms  
to construct test functions related to the topology of the surface. 

We point out that the similar ideas were used by Palmer \cite{Pa} to obtain a lower bound
for the index of the energy functional and also by Torralbo and Urbano in \cite{TU} to 
classify stable CMC spheres in homogeneous spaces. 

This paper is organized as follows. Section \ref{r} is devoted to state the main results of the paper. 
In Section \ref{preliminaries} we present the precise definition of index 
we discuss some known examples of CMC surfaces and they indices. 
In Section \ref{test} we present some auxiliares results that will be used in 
Section \ref{proofs}, a section dedicated to the proofs of the theorems. 

\section{Results}\label{r}

In this section we present the precise statements of our results.
For simplicity let us denote by $\bar{M}_c^{3}$ the space form of constant curvature $c\in \{0,1\}$,
that is,  the $3$-dimensional unit  sphere $\mathbb S^{3}$  for $c=1$ and the $3$-dimensional 
Euclidean space $\mathbb R^{3}$ for $c=0$.

Our main theorem read as follows. 

\begin{theorem}\label{ind}
Let $M^2$ be a compact constant mean curvature surface with genus $g$ immersed in $\bar M_c^{3}.$
Then, 
\[
\indw(M)\geq\frac{g}{3+c}.
\]
\end{theorem}	

Note that if $M$ is stable CMC surface in $\bar M_c^{3}$, then Theorem \ref{ind} implies
that $M$ is a topological sphere and by Chern-Hopf Theorem  \cite{C,H} it is a round sphere. 
So we recover Barbosa-do Carmo-Eschenburg Theorem \cite{BdCE} in  dimension $n=2$.

\begin{remark} It is a natural question to wonder whether there is a lower bound for the weak index 
of compact CMC hypersurfaces in $\mathbb S^{n+1}$ or in $\real^{n+1}$, $n\geq 3$, in terms of the first Betti number. 
\end{remark}

As a by product of the technique, we obtain a comparison between the eigenvalues 
$ \lambda_{\alpha}^J$ of the Jacobi operator of a  CMC  surface immersed in $\bar M_c^{3}$ 
and the eigenvalues $\lambda_{\beta}^{\Delta}$ of the Hodge  Laplacian acting on $1$-forms,
in the same spirit as Savo did for minimal hypersurfaces in $\mathbb S^{n+1}$. 

%
%
%
%
\begin{theorem}\label{esp}
Let $M^2$ be a compact immersed surface in $\bar M_c^{3}$ with constant mean curvature $H$.
Then for all positive integers $\alpha$ we have
\begin{eqnarray*} \lambda_{\alpha}^J\leq -2(c+H^2)+\lambda_{m(\alpha)}^{\Delta},
\end{eqnarray*}
where $m(\alpha)>2(3+c)(\alpha-1).$
\end{theorem}
 %
 %
 %

%
%
%

\begin{remark}
Bingqing Ma and Guangyue Huang obtained in \cite{MH}  a similar inequality as in Theorem \ref{esp}
for CMC hypersurfaces in $\mathbb S^{n+1}$, but involving the norm of the second fundamental form
of $M$. They comparison theorem does not imply index estimates for CMC hypersurfaces.
\end{remark}
%
%
%

\begin{remark} If $M$ is a compact minimal hypersurface in $\mathbb S^{n+1}$ which is not totally geodesic, then
$-n$ is an eigenvalue of $J$ with multiplicity at least $n+2$ (see \cite{U}, \cite{EL}, \cite{Sa}), in particular its index is at least
$n+3$, since the first eigenvalue is simple.
Savo \cite{Sa} used this fact to improve his estimates for the index of minimal hypersurfaces in the sphere.
We would like to point out here that Perdomo and Brasil proved in \cite{PB} that if $M$ is a compact CMC hypersurface in
$\mathbb S^{n+1}$ which is not totally umbilical, then $\indw(M)\geq n+1$, however the negative eigenvalues are not explicit and thus we cannot use it to improve our estimates. 
\end{remark}

\section{Preliminaries}\label{preliminaries}
In this section we will considerer hypersurfaces in any dimension. 

\subsection{The index of constant mean curvature hypersurfaces} 

Let $\bar M^{n+1}$ be a Riemannian manifold and
let  $\psi:M^n\rightarrow\bar M^{n+1}$ be an immersed two sided compact  hypersurface without boundary. 
We consider in $M$ the Riemannian metric $g$ induced
by  $\psi$. Let $\nabla$ and $\bar{\nabla}$ be  the Levi-Civita connections on $M$ and $\bar M$, respectively.
Fixed a unit normal vector field $N$ along  $M$, we will denote by $A$ its associated 
shape operator, that is,
\begin{equation*}
AX=-\bar{\nabla}_X N \ \ \mbox{for all} \ X\in TM.
\end{equation*}

The mean curvature function of $M$ is then defined as  $H=(1/n)tr A$.
It is well known that every smooth function $u\in C^{\infty}(M)$ induces a normal variation  
$\psi_t:M^n\rightarrow \bar{M}^{n+1}$
given by
\[
\psi_t(x)=exp_{\psi(x)}(tu(x)N_x),
\]
where $exp$ denotes the exponencial map in $\bar{M}^{n+1}.$
Since $M$ is closed and $\psi_0=\psi,$ there exists $\epsilon>0$ such that
\begin{equation*}
M_{u,t}=\{exp_{\psi(x)}(tu(x)N);x\in M\}
\end{equation*}
are immersed hypersurfaces  for all $t\in(-\epsilon,\epsilon).$
We can consider the area functional
$ \mathcal A_u:(-\epsilon,\epsilon)\rightarrow\mathbb{R} $ which is given by
\begin{equation*}
\mathcal A_u(t)=\int_MdM_{u,t},
\end{equation*}
where $dM_{u,t}$ is the $n$-dimensional area element of the metric induced on $M$ by $\psi_t.$
The first variation formula for the area is given by
\begin{equation*}
\mathcal A'_u(0)=-n\int_MuHdM.
\end{equation*}
As a direct consequence, minimal hypersurfaces are characterized as critical points of the area functional,
while constant mean curvature (CMC) hypersurfaces are the critical points of the area functional restricted to variations that preserves volume, that is,
$\displaystyle\int_Mu\,dM=0.$
For such critical points, the second variation of the area functional  is given by  the following 
quadratic form 
\[
\mathcal A''_u(0)=\int_M\|\nabla u\|^2 -(\overline {\ric}(N)+\|A\|^2)u^2 \,dM.
\]
Here $\|A\|^2=tr(A^2)$ is the Hilbert-Schmidt norm of $A$  and 
$\overline{\textrm{Ric}}(N)$ denotes the 
Ricci curvature of $\bar M$ in the direction of $N$. 
Integrating by parts we can write 
\[
\mathcal A''_u(0)=\int_MuJu \,dM.
\]
where $J=\Delta-\ric(N)-\|A\|^2$ is the so called Jacobi operator or stability operator of $M$.
We also point out that we are using the geometric definition of the Laplace-Beltrami operator, 
that is, $\Delta u=\delta du$
where $\delta w=-tr \nabla w$  for $w\in \Omega^1(M).$ 

The index of a CMC hypersurface $M$ is denoted by $\indw(M)$ and defined as the maximum dimension of any subspace $V$ of
\[
\displaystyle\mathcal{F}=\big\{u\in C^{\infty}(M); \int_Mu=0\big\}
\]
on which $\mathcal{A}''_u(0)$ is negative definite.
In other words, $\indw(M)$ is the number of negative eigenvalues of $J$, which is necessarily 
finite for closed hypersurfaces.	
In the special case that $\bar M^{3}$ has constant sectional curvature $c$, the Jacobi operator reads as 
\begin{equation*}\label{jeq}
J=\Delta-\|A\|^2-2c,
\end{equation*} and this form will be used in the rest of the paper.

%

\subsection{Examples of compact CMC hypersurfaces and their indices}
The most simple examples of closed CMC surfaces in the Euclidean space are 
the geodesic spheres, which are the only stables ones \cite{BdC}.
In 1986, Wente \cite{W} constructed the first examples of CMC tori in  $\real^3$ solving a question posed 
by Hopf in \cite{H}.
After that, all CMC tori in space forms were classified in a series of works 
(see \cite{Ab}, \cite{Bo}, \cite{PS} and \cite{UY}) 
and their indices were estimated by Lima, Sousa Neto and Rosmann in \cite{LSR} and by Rosmann in
\cite{Ro}. One can summarizing they results as follows 
\begin{quote}
{\it The index of a CMC torus in $\real^3$ is at least 8 and there are CMC tori with arbitrarily large index.}
\end{quote}
Many examples of compact CMC surfaces with genus $g\ge 2$ were constructed by Kapouleas in 
\cite{Ka1, Ka2, Ka3} using the gluing method.

When the ambient space is the round sphere $\mathbb S^{n+1}$ many examples of 
CMC hypersurfaces are known. Again, geodesic spheres are the only stable CMC hypersurfaces \cite{BdCE}. 
The next important family is given by the CMC Clifford tori, including the case $H=0$. 
 In fact, it was proved by S. Brendle  \cite{B} that that Clifford tori are the only minimal embedded tori in $\mathbb S^3$,  confirming longstanding conjecture of H. Lawson.
It is proved by B. Andrews and H. Li   \cite{AL}  that all  CMC embedded tori in $\mathbb S^3$  are rotational surfaces, confirming longstanding conjectures of  Pinkall-Sterling. We also note that Andrews and Li  \cite{AL}  gave a complete classifications of all  CMC embedded tori in 
$\mathbb S^3$.
It follows from the work of Simons \cite{S} that  any compact minimal hypersurface
not totally geodesic in $\mathbb S^{n+1}$ has index $\ind(M)\geq n+3$ and  it is also well known 
that the minimal Clifford torus has  index $n+3$. 
It is a natural problem to classify the minimal hypersurfaces $M\subset \mathbb S^{n+1}$ with $\ind(M)=n+3$. 
This problem was solved when $n=2$ by Urbano in \cite{U},  showing that the minimal Clifford tori are the only 
minimal surface in  $\mathbb S^3$ whose index is $5$. 
This problem is still open in higher dimensions. 

In the case $H\neq 0$, Perdomo and Brasil \cite{PB} proved that if $M\subset \mathbb S^{n+1}$ 
is a compact CMC hypersurface not totally umbilical, then $\indw(M)\geq n+1$ (see also \cite{A} for
a nice survey).
Recently, Alias and Piccione \cite{AP} showed the existence of infinite sequences of isometric 
embeddings  of tori with constant mean curvature in Euclidean spheres that are not isometrically 
congruent to the CMC Clifford tori, and accumulating at some CMC Clifford torus.
In the same spirit as above, the index of CMC tori of revolution in $\mathbb S^3$ 
were estimated by Rossman and
Sultana in \cite{RS} and by Ca\~nete in \cite{Ca}. On the other hand, 
higher genus CMC surfaces in $\mathbb S^3$ were constructed by 
Butscher-Pacard \cite{BP} (see also  \cite{B1,B2} for examples higher dimensions), but no index estimates 
is known.

\section{Test functions based on coordinates of vector fields} \label{test}
 
In this section we will consider  CMC immersed surfaces $\psi:M\to\bar{M}_c^{3}$.  
In the spherical case, $c=1$, we will  also consider  the unit normal vector field 
$\nu=-\psi$ along $\mathbb{S}^{3}$,  such that the second fundamental form 
of the inclusion map, $\bar{\psi}:\mathbb{S}^{3}\rightarrow \mathbb{R}^{4}$, is the identity.
That is,  $D_X\nu=-X$ where $D$  denotes  the Levi-Civita connection in the Euclidean space. 
It follows immediately that
 \begin{eqnarray*}\label{conexao}
D_Y X-\nabla_Y X=\left\langle X,Y\right\rangle\nu+\left\langle AX,Y\right\rangle N,  
 \mbox{for all} \  X,Y\in TM
\end{eqnarray*}
for $M$ immersed in $\mathbb{S}^{3}$ and
 \begin{eqnarray*}\label{conexao_r}
D_Y X-\nabla_Y X=\left\langle AX,Y\right\rangle N,  \ \mbox{for all} \  X,Y\in TM
\end{eqnarray*}
for $M$ immersed in $\mathbb{R}^{3}.$

Fixed an orthonormal basis  $\mathcal{E}=\{\bar E_1,\dots,\bar E_{3+c}\}$  of parallel vector fields 
on $\mathbb{R}^{3+c}$ we will denote by 
\[{E_i}:=\bar E_i-\z \bar E_i,N\x N-c\z \bar E_i,\nu\x\nu 
\]
the vector fields given by the  orthogonal projection of $\bar E_i$ on $TM$. 
For the fixed basis $\mathcal{E}$ we will consider the
smooth support functions $f_i, g_i:M\to\real$ given by
\[
f_i=\left\langle \bar{E_i},\nu\right\rangle \quad  \textrm{ and }  \quad g_i=\left\langle \bar{E_i},N\right\rangle,
\]
for $1\leq i\leq 3+c$.

Let $\xi\in TM$ be a smooth vector field on $M$ and let $\omega$ denote its dual $1$-form,
that is $\xi = \omega^\#$.
Inspired in the works of Ros and Savo we will use the coordinates of $\xi\in TM$
as test functions. They are given by
\[
w_i:=\z  E_i,\xi\x, \quad 1\leq i\leq 3+c.
\] 

Let us denote by $\nabla^{*}\nabla$  the rough Laplacian acting on vector fields and by $\Delta$ 
the Hodge Laplacian acting on $1$-forms. They are defined, respectively, by
\[
\nabla^{*}\nabla\xi=-\trace \nabla^2\xi \quad \mbox{and}\quad \Delta\omega=d\delta\omega+\delta d\omega,
\]
where  $d$ is the exterior differential and  $\delta=-\star d\star$ is the formal adjoint of $d$ with respect to the canonical $L^2$-inner product on $1$-forms induced by the Riemannian metric of $M$.
We define  the Laplacian of the vector fields $\xi$ as being $\Delta\xi = (\Delta\omega)^{\#}$ and so
these Laplacians are  related by the well known Bochner formula 
\begin{equation}\label{bochner}
\Delta \xi=\nabla^{*}\nabla\xi+K\xi,
\end{equation}
where $K$ is the Gauss curvature of $M$.

In order to compute the Jacobi operator of $w_i$ we need the following lemma.
\begin{lemma}\label{lapwi}Let $M^{2}$ be an orientable CMC  surface immersed in $\bar M_c^3.$ 
Then, using the above notation we have 
\[
\Delta w_i=(\|A\|^2-4H^2)w_i+2H\z AE_i,\xi\x-2g_i\z A,\nabla \xi\x+2cf_i\diver\xi+\z E_i,\Delta \xi\x,
\]
for $1\leq i\leq 3+c$. 
\end{lemma}
\proof  Fixed a point $p\in M$, we consider a local orthonormal frame $\{e_1,e_2\}$ on $M$ which is geodesic at $p.$ 
A direct computation shows that  $\nabla_{e_\ell}E_i=g_iAe_\ell+cf_ie_\ell,$ 
$\ell=1,2.$
Thus, using Einstein summation notation, we get
\begin{eqnarray*} 
\Delta w_i &=&-e_\ell e_\ell\z  E_i,\xi\x\\
                  &=&-e_\ell(\z \nabla _{e_\ell} E_i,\xi\x+\z  E_i, \nabla_{e_\ell}\xi\x)\\
                 &=&-\z\nabla_{e_\ell}\nabla_{e_\ell}E_i,\xi\x-2\z \nabla_{e_\ell} E_i,\nabla_{e_\ell}\xi\x-
                   \z E_i, \nabla_{e_\ell}\nabla_{e_\ell}\xi\x\\
                   &=&-\z \nabla_{e_\ell}(g_iAe_\ell+cf_ie_\ell),\xi\x-2\z g_iAe_\ell +
                   cf_ie_\ell,\nabla_{e_\ell}\xi\x+\z E_i,\nabla^*\nabla \xi\x\\
                 &=&\z AE_i,e_\ell\x\z Ae_\ell,\xi\x-g_i\z (\nabla_{e_\ell}A)e_\ell,\xi\x+c\z E_i,e_\ell\x\z e_\ell,\xi\x\\
                 &&-2g_i\z A,\nabla \xi\x+2cf_i\diver \xi+\z E_i,\nabla^*\nabla \xi\x\\  
                 &=&\z A^2E_i,\xi\x-2g_i\z A,\nabla \xi\x+2cf_i\diver\xi
                 +\z E_i,\Delta \xi\x\\&&-2H^2\z E_i,\xi\x+\frac{\|A\|^2}{2}\z E_i,\xi\x,
                      \end{eqnarray*}
where in the last equality we used the Bochner equation (\ref{bochner}) and the Gauss equation in the form
below:
\[
K=c+2H^2-\frac{1}{2}\|A\|^2.
\]
To conclude the proof we note that shape operator $A$ satisfies the following equation
\[
A^2=\frac{1}{2}(\|A\|^2-4H^2)I_2+2HA.
\] 
\endproof

We end this section noting that the coordinates of harmonic vector fields are 
admissible functions to compute the index of CMC surfaces. Moreover we have:

\begin{lemma}\label{wi} If $\xi\in TM$ is a harmonic vector field, then 
$w_i= \langle E_i, \xi\rangle$ and  $\bar w_i=\langle E_i,\star\xi\rangle$ satisfy
\[
\int_M w_i =\int_M \bar w_i =0
\]
for $1\leq i\leq 3+c.$
\end{lemma}
\proof
If  $\xi$ is harmonic we have $\diver\xi=0$ and therefore
\[		
\int_M w_i=-\int_M\langle \nabla f_i,\xi\rangle=-\int_Mf_i\diver\xi=0.
\]
Now, recall that, in an orthonormal basis $\{e_1,e_2\}$ of $TM$, the Hodge star operator is defined by 
\[
\star e_1=e_2,\star e_2=-e_1.
\]
Since $\Delta$ comutes with $\star$ it follows that $\star\xi$ is also a harmonic vector field on $M$
and this concludes the proof.
\endproof


\section{Proofs of the theorems}\label{proofs}

For simplicity we will present the proofs in the case of CMC surfaces in the unit sphere $\mathbb S^3$. 
The case of CMC surfaces in  $\mathbb R^3$ follows the same steps.

Let us denote by   $\mathcal{L}^{\Delta}_m$ the vector space given by the direct sum 
of the  eigenspaces generated by $\xi_1, \xi_2, \ldots, \xi_m$, 
the first $m$ eigenfunctions of the Hodge Laplacian $\Delta$ and
let us denote by $\mathcal{H}^1(M)$  the vector space of the harmonic vector fields on $M$. 
Notice that $\dim \mathcal{H}^1(M) =2g$, where $g$ is the genus of surface $M$, 
and $\mathcal{H}^1(M)\subset \mathcal{L}^{\Delta}_m$ .
%

\subsection{Proof of Theorem \ref{esp}}

Since $J$ is an elliptic self-adjoint operator, it admits a sequence of eigenvalues diverging to infinity,
\[
\lambda_1^{J}\leq \lambda_2^{J}  \leq\cdots\leq \lambda_k^{J} \leq\cdots
\]
Fix an orthonormal basis  $\{\phi_1,\phi_2,\ldots\}$ of  $C^{\infty}(M)$ 
given by eigenfunctions of the Jacobi operator, that is, $J\phi_i=\lambda_i^J\phi_i$.
We denote by  $\mathcal{J}^{p} :=\z \phi_1,\cdots,\phi_{p}\x^{\bot}$
the linear space orthogonal to the first $p$  eigenfunctions of the  Jacobi operator.

Initially, we look for vector fields $\xi\in \mathcal{L}^{\Delta}_m$  such that the functions 
$w_i,\bar w_i\in \mathcal{J}^{\alpha-1}$, for some $\alpha\in \mathbb N$ and $i\in\{1,\dots, 4\}.$ In other words, we have a system with $8(\alpha-1)$  homogenous linear equations
in the variable $\xi$
\begin{equation}\label{sys}
\int_Mw_i\phi_k=\int_M\bar w_i\phi_k=0,\quad 1\leq i\leq 4 \quad \mbox{and} \quad 1\leq k\leq \alpha-1.
\end{equation}
Therefore, if $m(\alpha)=\dim \ \mathcal{L}^{\Delta}_m>8(\alpha-1),$ then the system (\ref{sys})  has at least a  non trivial solution $\xi\in \mathcal{L}^{\Delta}_m$  such that $w_i,\bar w_i\in \mathcal{J}^{\alpha-1}$ for all $1\leq i\leq 4.$ By the min-max principle we have
\begin{eqnarray*}
\lambda_{\alpha}^J\leq\frac{\int_Mw_iJw_i}{\int_Mw_i^2}\quad \mbox{and}\quad \lambda_{\alpha}^J\leq\frac{\int_M\bar w_iJ\bar w_i}{\int_M\bar w_i^2}. 
\end{eqnarray*}
Now, using Lemma \ref{lapwi} we get 
{\setlength\arraycolsep{1pt}
\begin{eqnarray*}
\lambda_{\alpha}^J\int_Mw_i^2&\leq &-(2+4H^2)\int_Mw_i^2+2H\int_M\z E_i,A\xi\x w_i\\
&& +\int_M\z E_i,\Delta \xi\x w_i-2\int_Mg_i\z A,\nabla \xi\x w_i+2\int_Mf_i\delta\xi w_i.
\end{eqnarray*}}

Summing upon $i=1,\ldots, 4$ we obtain 
{\setlength\arraycolsep{1pt}
\begin{eqnarray*}
\lambda_{\alpha}^J\int_M\|\xi\|^2&\leq &-(2+4H^2)\int_M\|\xi\|^2+2H\int_M\z A\xi,\xi\x +\int_M\z\Delta \xi,\xi\x.
\end{eqnarray*}}

Analogously, we do the same to the test functions $\bar w_i$:
{\setlength\arraycolsep{1pt}
 \begin{eqnarray*}
\lambda_{\alpha}^J\int_M\|\xi\|^2&\leq &-(2+4H^2)\int_M\|\xi\|^2+2H\int_M\z A\star\xi,\star\xi\x +\int_M\z\Delta \star\xi,\star\xi\x.
\end{eqnarray*}}

Summing these last two inequalities we have 
{\setlength\arraycolsep{0pt}
 \begin{eqnarray}\label{lambda}
\lambda_{\alpha}^J\int_M\|\xi\|^2&\leq &-(2+4H^2)\int_M\|\xi\|^2+H\int_M\z A\xi,\xi\x+\z A\star\xi,\star\xi\x\\
						&&+\frac{1}{2}\int_M(\z\Delta \xi,\xi\x+\z\Delta \star\xi,\star\xi\x) \nonumber.
    \end{eqnarray}}
Now, we observe that 
\begin{equation}\label{AA}
\z A\xi,\xi\x+\z A\star\xi,\star\xi\x=2H\|\xi\|^2
\end{equation} 
for any $\xi\in TM$.   
If $\xi\in \mathcal{L}^{\Delta}_m$ we can  $\xi=\alpha_i\xi_i$ and then
 \begin{eqnarray}\label{xi}
\int_M\z\Delta \star\xi,\star\xi\x&=&\int_M\z\Delta \xi,\xi\x\\ \nonumber
                                               &=&\lambda_i^{\Delta}\int_M\alpha_i\alpha_k \z\xi_i,\xi_k\x\\ \nonumber
                                               &\leq&\lambda_{m(\alpha)}\int_M\|\xi\|^2. \nonumber
    \end{eqnarray}
Plugging (\ref{AA}) and (\ref{xi}) into (\ref{lambda}) we obtain
\[
\lambda_{\alpha}^J\leq -2(1+H^2)+\lambda_{m(\alpha)}^{\Delta}
\]
where $m(\alpha)>8(\alpha-1)$.

\subsection{Proof of Theorem \ref{ind}}
We start as in proof of Theorem \ref{esp} but now choosing an orthonormal basis  $\{\phi_1,\phi_2,\cdots\}$  
of the  space  $\mathcal{F}=\big\{u\in C^{\infty}(M); \int_Mu=0\big\}$ 
given by eigenfunctions of the Jacobi operator. 
We already know from Lemma \ref{wi} that for any $\xi\in\mathcal{H}^1(M)$ the test functions 
$w_i,\bar w_i\in\mathcal{F}.$  

Consider then vector fields $\xi\in \mathcal{H}^1(M)$  such that the test functions $w_i,\bar w_i\in \mathcal{J}^{\alpha-1}$, for some $\alpha\in \mathbb N$ and $i\in\{1,\dots, 4\}.$ 
As before,   if $\dim \mathcal{H}^1(M)=2g>8(\alpha-1),$ then the system (\ref{sys})  has at least a 
non trivial solution $\xi\in \mathcal{H}^1(M)$.
Following the same steps as above, we use Lemma \ref{lapwi} but now for the harmonic vector fields $\xi$
and its dual $\star \xi$. We  obtain 
%
%
%
{\setlength\arraycolsep{1pt}
 \begin{eqnarray*}
\lambda_{\alpha}^J\int_M\|\xi\|^2&\leq &-2(1+H^2)\int_M\|\xi\|^2.
    \end{eqnarray*}}

Hence,  we conclude that $\lambda_{\alpha}^J< 0$ and then $\indw (M)\geq \alpha.$ 
Since $\alpha$ can be chosen as the largest integer such that $2g>8(\alpha-1)$ we get 
\[
\indw(M)\geq \frac{g}{4}.
\]


\bibliographystyle{amsplain}

\end{document}